\documentclass[11pt]{amsart}
\usepackage{amsfonts, amstext, amsmath, amsthm, amscd, amssymb}
\usepackage{epsfig, graphics, psfrag}
\newtheorem{theo}{Theorem}[section]
\newtheorem{pro}[theo]{Proposition}
\newtheorem{lem}[theo]{Lemma}
\newtheorem{defi}[theo]{Definition}

\newtheorem{remark}[theo]{Remark}
\newtheorem{corol}[theo]{Corollary}
\newtheorem{exam}[theo]{Example}

\def \proof {{\bf Proof$\colon$}\ }

\def \Z{{\bf Z}}

\begin{document}

\title[]{Cosmetic crossing changes of fibered knots}
\author[E. Kalfagianni]{Efstratia Kalfagianni}

\address[]{Department of Mathematics, Michigan State
University,
East Lansing, MI, 48824}

\email[]{kalfagia@math.msu.edu}
%\thanks{ {\today}}
\thanks{ {Research partially supported by the NSF}}
%{{NSF FRG/DMS-0456155,
%DMS-0805942 and by a grant from the 
%IAS}}

\thanks{Authors e-mail address: kalfagia@math.msu.edu}
\thanks{ \today}
\maketitle

\begin{abstract} We prove the nugatory crossing conjecture for fibered knots.
 We also show that if a knot $K$ is 
$n$-adjacent to a fibered knot $K'$, for some $n>1$, then either the genus of 
$K$ is larger than that of $K'$ or $K$ is
isotopic to $K'$.

\smallskip
\smallskip

\smallskip
\smallskip

\smallskip
\smallskip

\noindent {\it Keywords:} crossing change, commutator length of a Dehn twist, fibered knot, mapping class group, Heegaard
splitting, Thurston norm.
\smallskip
\smallskip

\smallskip
\smallskip

\noindent {\it Mathematics Subject Classification:} 57M25, 57M27, 57N65, 57R30, 20F65.
\end{abstract}

\smallskip
\smallskip

%\tableofcontents

\section{Introduction}

\smallskip

An open question in classical knot theory is the question of when a crossing 
change
on a knot changes the isotopy class of the knot. The purpose of this paper 
is to answer this question for fibered knots.

A crossing disc for a knot $K\subset S^3$
is an embedded disc $D\subset S^3$
such
that $K$ intersects ${\rm int}(D)$ twice with
zero algebraic intersection number. A crossing change on $K$ can be achieved
by twisting $D$ or equivalently by performing appropriate Dehn surgery of $S^3$
along the crossing circle $\partial D$.
The crossing is called nugatory if and only if
$\partial D$ bounds an embedded disc in the complement of $K$.
This disc and $D$ form a
2-sphere that decomposes $K$ into a connected sum,
where some of the summands may be trivial.
Clearly, changing a nugatory crossing doesn't
change the isotopy class of a knot.
Problem 1.58  of \cite{kn:kirby} asks whether the converse is true (see also 
\cite{kn:to} for related conjectures):
That is,  if a crossing change on a knot $K$ yields a knot
isotopic to $K$ is the crossing
nugatory? 

In the case that $K$ is the trivial knot an affirmative answer follows from work
of Gabai \cite{kn:ga}. 
An  affirmative answer is also known
in the case of 2-bridge knots \cite{kn:to}.
In this paper we will show the following.

\begin{theo} \label{theo:nugatory} Let $K$ be a fibered knot.
A crossing change on $K$ yields a knot
isotopic to $K$ if and only if the crossing is
nugatory.
\end{theo}

To give a brief outline
of the proof, let $K$ be a  
fibered knot such that a crossing 
change on $K$ gives a knot $K'$ that is isotopic to $K$.
The complement of $K$ is fibered over $S^1$ with fiber, say $S$; a minimal genus Seifert
surface of $K$. A result of Gabai implies that the crossing change from $K$ to $K'$
can be achieved along an arc that is properly embedded
on $S$. Equivalently, the crossing change can be achieved by twisting
$K$ along
a meridian disc $D$ of a handlebody
neighborhood $N$
of the fiber. 
Using geometric properties of fibered knot complements,
the problem reduces to the question of whether a power of a Dehn twist on 
the surface
$\partial N$ along
he curve $\partial D$, can be written as a single commutator in the 
mapping class group of the surface.
A result of
Kotschick implies that a product
of Dehn twists of the same sign, along a collection of disjoint, homotopically essential curves on an orientable surface 
cannot be written as a
single commutator in the mapping class group of the surface. 
Using this result,
we show that the assumption that $K$ is isotopic
to $K'$ implies that
$\partial D$ bounds a disc in the complement of $K$.
\vskip 0.04in

Theorem \ref{theo:nugatory} says that an essential crossing change always 
changes the isotopy class of a fibered knot. It is natural to ask whether
the crossing change produces a simpler or
more complicated knot with respect to some knot complexity. 
A complexity function whose interplay with crossing changes has been 
studied using the theory of taut foliations and sutured 3-manifolds 
is the knot genus. Simple examples show that a 
single crossing change
may decrease or increase  the genus of a knot even 
if one stays within the class of fibered knots. However there are interesting consequences
if one replaces 
a crossing change by the more refined notion of knot adjacency 
\cite{kn:kl}, \cite{kn:kl1} :
We recall that $K$ is called 
2-adjacent to $K'$
if $K$ admits a projection that contains two crossings such that changing 
any of them or both of them simultaneously, transforms $K$ to $K'$.
\begin{theo}\label{theo:unique1}
Suppose that $K'$ is a fibered knot and that $K$ is 2-adjacent to $K'$. Then either $K$ is isotopic to $K'$ or $K$ has a strictly larger genus than $K'$.
\end{theo}
%It is  known that the Alexander polynomial can be used to
%detect non-fibered knots. Namely, the Alexander polynomial of a fibered knot is monic.
%Theorem \ref{theo:unique1}
%can be used to obtain
%criteria for detecting non-fibered knots when the Alexander polynomial
%provides inconclusive evidence. This direction is explored in \cite{kn:kl1}
%where we also provide
%applications in the theory of finite type invariants of 3-manifolds.

We organize the paper as follows:
In Section 2 we summarize the mapping class group results that we need for
the proof of Theorem \ref{theo:nugatory} and in Section 3 
we summarize known properties of fibered knot complements.
In Section 4
we discuss a setting relating fibrations  of knot
complements and Heegaard splittings of $S^3$,
from the point of view needed in the rest of the paper.
In Section 5, we study nugatory crossings of fibered knots and we prove
Theorem \ref{theo:nugatory}.
In Section 6 we study adjacency to fibered knots and
prove  Theorem \ref{theo:unique1}.

Throughout the paper we work in the PL
or the smooth category.

\smallskip

\section{Commutator length  and Dehn twists}

\subsection{Commutators in the mapping class group.} 
Let $\Sigma_k$ denote a closed oriented surface of genus $k$ and let
$\Gamma_k$ denote the mapping class group of $\Sigma_k$. That is $\Gamma_k$
is the group of isotopy classes  of orientation preserving  homeomorphisms 
$\Sigma_k \longrightarrow \Sigma_k$.
Let $\Gamma'_k:=[\Gamma_k, \ \Gamma_k]$ denote the commutator subgroup of $\Gamma_k$.
An element  $f\in \Gamma'_k$ is written as a
product of commutators.
The commutator length of $f$, denoted by
$c(f)$, is the minimum number of factors needed to express $f$ as a product
of commutators. In the recent years, the growth of the commutator length of Dehn twists
has been studied using methods from the theory of symplectic four-manifolds \cite{kn:eko},
\cite{kn:brko}, \cite{kn:ko}, \cite{kn:kOz}.
In this paper, we will need a result of D. Kotschick which we recall below.

For a simple closed curve $a\subset \Sigma_k$ let $T_a$ denote the right
hand Dehn twist about $a$; then the left hand Dehn twist about $a$ is $T_a^{-1}$. 

\begin{theo} \label{theo:length} {\rm [Theorem 7, \cite{kn:ko}] }
Let $\Gamma_k$ be the mapping class group of a closed oriented surface $\Sigma_k$ of
genus $k\geq 2$. Suppose that $a_1,\ldots, a_m \subset \Sigma_k$ are
homotopically essential, disjoint, simple closed
curves on $\Sigma_k$. Let $f:= T_{a_1} \cdot T_{a_2}\cdot \ldots \cdot T_{a_m}$ denote the product
of
right-handed Dehn twists along $a_1,\ldots, a_m$. Suppose that for some $q>0$ we have   $f^q=T_{a_1}^{q} \cdot T^{q}_{a_2}\cdot \ldots \cdot T^{q}_{a_m}  \in \Gamma'_k$.
Then, we have
$$ c(f^q)\geq 1+\frac{qm}{18k-6}.$$
\end{theo}
\smallskip

We will need the following corollary of Theorem \ref{theo:length}.

\begin{corol} \label{corol:one} Let $\Gamma_k$ be the mapping class group of
a closed oriented surface $\Sigma_k$ of genus $k\geq 2$. Let $a\subset 
\Sigma_k$ be
a simple closed curve.
Suppose that there exist $g,h \in \Gamma_k$  such that
$$T_a^q= [g,h]=ghg^{-1}h^{-1},$$
for some $q\neq 0$. Then $a$ is homotopically trivial on $\Sigma_k$.
\end{corol}
\smallskip

The proof of Theorem \ref{theo:length} given in \cite{kn:ko} relies on the
theory of Lefschetz fibrations,
which, as the author points out, is sensitive to the chirality of Dehn twists.
In fact,
the argument
of \cite{kn:ko} breaks down if one allows $f$ to be a product of
right-handed Dehn twists and their inverses and,
as the following example shows, 
Theorem \ref{theo:length} is not true in this case.
In subsequent sections we will discuss how this situation is reflected when 
one tries to apply Theorem \ref{theo:length} to the study of crossing changes
that do not alter the isotopy class of fibered knots (see Example \ref{exam:figure8}).

\begin{exam} \label{exam:more} {\rm [Example 9, \cite{kn:ko}]}
{\rm Suppose that $a\subset \Sigma_k$ is an essential simple closed loop on
a closed oriented surface of genus at least two. Let $g: \Sigma_k \longrightarrow \Sigma_k$  be
an orientation preserving homeomorphism such that $a\cap g(a)=\emptyset$.
We will also use $g$ to denote the mapping class of $g$. Set $b:=g(a)$ and
set $f:=T_a T_b^{-1}$.
In the mapping class group $\Gamma_k$ we have $gT_ag^{-1}=T_{g(a)}$
or equivalently  $T_a=g^{-1}T_bg$.
Since $a,b$ are disjoint we also have $T_aT_b^{-1}=T_b^{-1}T_a$. Thus
$$f^{q}=(T_aT_b^{-1})^q=T_a^qT_b^{-q}=(g^{-1}T_bg)^q T_b^{-q}=[g^{-1}, T_b^{q}],$$
for all $q>0$. Hence we have  $c(f^q)=1$ showing that
Theorem \ref{theo:length} is  not true in this case.}
\end{exam}

\subsection{When is $T_a^q$ trivial?}
It is known that if $a$ is homotopically essential on $\Sigma$
then no non-trivial power $T_a^q$ ($0 \neq q\in \Z)$ is isotopic to the identity on
$\Sigma_k$. This statement is well known to researchers working on mapping class groups: It is for example
asserted in  \cite{kn:blm} when the authors state that the kernel of the reduction homomorphism corresponding
to an essential simple closed curve $a$ is the free abelian group generated by Dehn twists along $a$. Below we include a
proof that
uses properties of intersection numbers stemming from Thurston's study of surface homeomorphisms \cite{kn:trav}.

\begin{pro} \label{pro:nevertrivial} Suppose that $T_a^q=1$ in the mapping class group $\Gamma:=\Gamma(\Sigma_k)$, $k>0$. Then,
either $a$ is homotopically trivial on $\Sigma_k$ or $q=0$.
\end{pro}
\proof Suppose that the curve $a$ is not homotopically trivial on $\Sigma$ and that $T_a^q=1$ in the mapping class group
$\Gamma(\Sigma)$. We will argue that $q=0$.
First suppose that $a$ is a non-separating loop on $\Sigma$. Then, we can find an embedded loop $b$ that intersects $a$ 
exactly once. Orient $a$, $b$ so that the algebraic intersection of $a,b$ is 1; that is $<a,b>=1$.
In $H_1(\Sigma)$ we have
$$T_a^q(b)=b+q<a,b>a=b+qa.$$
\noindent  Thus  we have
$<T_a^q(b), b>= <b, b>+q<a,b>$ which, since $T_a^q(b)=b$, gives $q=0$ as desired.
If $a$ is separating, we  appeal to the geometric intersection number. For $b$ a simple closed loop
on $\Sigma$ let $i(a,b)$ denote the intersection number; the minimum number of intersections in the isotopy classes of $a$ and $b$.
Since we assumed that $a$ is homotopically essential on $\Sigma$, we can find $b$ so that $i(a,b)\neq 0$.
By Expos\'e 4 \cite{kn:trav} we have the following:

$$i(T_a^q(b),\  b)= |q|\  (i(a, \ b))^2.$$
Since $T_a^q=1$, we have $0=i(b,\  b)=i(T_a^q(b),\  b)$.
Thus we obtain
$|q|\  (i(a, \ b))^2=0$; which implies that $q=0$. \qed
\smallskip

{\bf Notation.} To simplify our notation, throughout the paper, we will use $\Sigma:=\Sigma_k$ to denote an oriented
surface
of any genus $k\geq 0$ and $\Gamma:=\Gamma_k$ to denote the mapping class group of $\Sigma$.
Also, as we've done in this section, we will use the same symbol to denote a homeomorphism
of $\Sigma$ and its class in $\Gamma$.
\smallskip

\section{Uniqueness properties of knot fibrations}
Here we summarize some known properties of fibered knots that we need in subsequent sections.
For details and proofs the reader is referred to Section 5 of \cite{kn:bz} and \cite{kn:Waldhausen}.
Suppose that $K$ is a fibered knot and let $S$ be a minimum
genus Seifert surface
for $K$. Let ${\eta(K)}$ denote a tubular neighborhood of $K$. Then
the complement $\overline {S^3\setminus {\eta(K)}}$ admits a fibration
over $S^1$ with fiber $S$. More specifically, it is shown that the complement
$\overline {S^3\setminus {\eta(K)}}$  cut along $S$ is homeomorphic to $S\times [-1 , 1]$. Thus,
there is an orientation preserving homeomorphism $h: S \longrightarrow S$
such that $\overline {S^3\setminus {\eta(K)}}$  is obtained from $S\times [-1 , 1]$ by identifying
$S\times \{ -1\}$ with $S\times \{ 1\}$ so that $(x,-1)=(h(x), 1)$. 
The map $h$ is called the monodromy of the fibration.
We write $$\overline {S^3\setminus {\eta(K)}}= {S \times J}/ {h},$$
\noindent where $J:=[-1, 1]$.  We need the following:

\begin{pro} \label{pro:510} a) Let $M:=\overline {S^3\setminus {\eta(K)}}= {S \times J}/ {h}$ be an oriented, fibered knot complement
and set $S_1:=S\times \{1\}=S\times \{-1\}$.
Given a  minimum genus Seifert surface
$S_2$, with $\partial S_2=\partial S_1$, there 
there exists an orientation preserving homeomorphism of $M$ that is fixed on $\partial M$
and brings
$S_2$ to the fiber $S_1$. In fact such a homeomorphism is isotopic to the identity on $M$ by an isotopy relatively the boundary $\partial M$.

b) Let $M:={{S \times J}/ {h}}$
and $M':={{S' \times J}/ {h'}}$ be
fibered,  oriented knot complements.
Then, there exists
an orientation preserving
homeomorphism $F: M \longrightarrow M'$, with $F(\partial S\times \{j\})=\partial S'\times \{j\}$ ($j\in J$)
if and only if there exists an orientation preserving surface homeomorphism
$f: (S, \partial S) \longrightarrow (S', \partial S')$ such that
$fhf^{-1}$ and $h'$ are equal up to isotopy on $S'$.
\end{pro}

\smallskip

\section{Splittings of fibered knot complements}
Given a fibration of a knot complement $M:=\overline {S^3\setminus {\eta(K)}}={S \times J}/ {h}$,
set $N_1:= S \times [0, 1]$, $N_2:=
S \times [-1, 0]$,
$E:= \partial S \times (0, 1)$ and $E':= \partial S \times (-1, 0)$. We
have $\partial N_1= (S\times \{0\})\cup E \cup (S\times \{1\})$. Similarly,
we have $\partial N_2= (S\times \{-1\})\cup E' \cup (S\times \{0\})$.
We will assume that $K:= \partial S\times \{{{1}\over {2}}\}$ on
$\partial N_1$.
Define 
$g: \partial N_1\longrightarrow \partial N_1$ by
$$g(x, 0)=(x,0), \ \ {\rm for} \ \ x\in S, \eqno (4.1) $$
$$g(x, t)=(x, t), \ \ {\rm for} \ \ x\in \partial S \ \ {\rm and}\ \ 0<t<1, \eqno (4.2) $$
$$g(x, 1)=(h(x),1), \ \ {\rm for} \ \ x\in S. \eqno (4.3)$$
Consider the homeomorphism $rg: \partial N_1\longrightarrow \partial N_2$, where
$r: N_1 \longrightarrow N_2$ is defined  by $(x,t)\rightarrow (x, -t)$. 
We obtain a Heegaard splitting 
$$
S^3=N_1{ \cup}_{rg} N_2:=N_1{ \sqcup} N_2 / \{ y\sim rg(y) \ | \  y\in \partial N_1\}, \eqno (4.4)$$
such that
$K$ lies on the Heegaard surface.  Next we push $K$ on $S\times \{{{1}\over {2}}\}$ slightly in the interior of $N_1$
and then we take  $A(K)$ to be an annulus neighborhood of $K$ on $S\times \{{{1}\over {2}}\}$.
Next we remove a
tubular neighborhood of $K$, say $\eta(K):= A(K)\times (\{{{1}\over {2}}\}- \epsilon, \  \{{{1}\over {2}}\}+\epsilon)$,
from
${\rm int}(N_1)$ and we set $H_1:=\overline {N_1\setminus {\eta(K)}}$. 
The decomposition
$$
M=H_1{ \cup}_{rg} N_2:=H_1{ \sqcup} N_2 / \{ y\sim rg(y) \ | \  y\in \partial N_1\}, \eqno (4.5)$$
is called
the $HN$-splitting corresponding to the fibration of $M$.
The  $HN$-surface
of this decomposition is
$Q:=\partial N_1{ \sqcup} \partial N_2 / \{ y\sim rg_1(y) \ | \  y\in \partial N_1\}$.
Now set $N:=N_1=S\times [0,1]$ and 
identify $N_2$ with $(-N)$  via  $r^{-1}$, where  $(-N)$ denotes $N$ with the opposite orientation. Also set $H:=\overline {N\setminus {\eta(K)}}$ 
and  $\Sigma:=\partial N_1$ and let $i: N\longrightarrow (-N)$ denote the orientation reversing involution. 

\begin{defi} \label{defi:modelfib}
The pair $(\Sigma,\   g)$ is called the $HN$-model
associated to the fibration $M={{S \times J}/ {h}}$. Note that, by (4.1)-(4.4), $g$
is the identity on $\Sigma\setminus (S\times \{1\})$.
\end{defi}

\begin{defi} \label{defi:model} Let $K$ be fibered knot  with $M:=\overline {S^3\setminus {\eta(K)}}={S \times J}/ {h}$
and let $H, N$ and $\Sigma$ be as above.  Also let $g_1: \Sigma \longrightarrow \Sigma$ be an orientation preserving
homeomorphism. The pair $(\Sigma,\ g)$ is called an $HN$-model  for $M$ if there is an orientation-preserving homeomorphism
$\Phi: M \longrightarrow H{ \cup}_{g_1} (-N)$, such that $\Phi| \partial \eta(K)= {\rm id}$.
Here, $$H{ \cup}_{g_1} (-N):=H{ \sqcup} (N) / \{ y\sim ig_1(y) \ | \  y\in \Sigma\}.$$
The surface $\Sigma{ \sqcup} \Sigma / \{ y\sim ig_1(y) \ | \  y\in \Sigma\}$ will be called the $HN$-surface
of the decomposition $H{ \cup}_{g_1} (-N)$
\end{defi}

The next Lemma reformulates Part(b) of Proposition \ref{pro:510} in terms of the models of 
the two fibrations.

\begin{lem} \label{lem:model510}  Let $M:={{S \times J}/ {h}}$
and $M':={{S' \times J}/ {h'}}$ be
fibered,  oriented knot complements in $S^3$ and  let $(\Sigma,\   g)$, $(\Sigma,\   g')$ denote the  models corresponding to the fibration of $M$, $M'$ respectively. There exists an orientation-preserving
homeomorphism $F: M \longrightarrow M'$, with $F(\partial S\times \{j\})=\partial S'\times \{j\}$ ($j\in J$)
if and only if
there is an orientation-preserving homeomorphism $f: \Sigma \longrightarrow \Sigma$
such that in the mapping class group $\Gamma=\Gamma(\Sigma)$ we have
$$g'=fg f^{-1}. $$
\end{lem}
\proof By Proposition \ref{pro:510} there exists
an orientation-preserving
homeomorphism $F: M \longrightarrow M'$, with $F(\partial S\times \{j\})=\partial S'\times \{j\}$ ($j\in J$)
if and only if there exists an orientation-preserving surface homeomorphism
$f: (S, \partial S) \longrightarrow (S', \partial S')$ such that
$fhf^{-1}$ and $h'$ are equal up to isotopy on $S'$. Now $g$ is constructed out of $h$ as in (4.1)-(4.4); in a similar fashion $g'$ is constructed out of $h'$. Set $I:=[0, \ 1]$.
We may extend
$f$ to a homeomorphism of pairs $(S\times I, \partial(S\times I))\longrightarrow (S'\times I, \partial(S'\times I))$ by defining $f(x, t)=(f(x), t)$. By our construction of 
the $HN$-splittings corresponding to fibrations this extension is considered as a map $(N, \Sigma) \longrightarrow (N, \Sigma)$.
Since, $g$ is the identity on $\Sigma\setminus (S\times \{1\})$ and 
$g'$ is the identity on $\Sigma\setminus (S'\times \{1\})$, we have 
$g'=fg f^{-1}$ up to isotopy on $\Sigma$. \qed

\smallskip

\smallskip

Let $Q$ denote the   $HN$-surface of the
splitting associated to a fibration $\overline {S^3\setminus {\eta(K)}}={S \times J}/ {h}$. By construction we have a surface $S_1\subset S\times \{{{1\over 2}}\}$
that is disjoint from $Q$. Furthermore $S_1$ and $S\times \{{{1\over 2}}\}$ differ by an annulus.
We will think of this $HN$-surface as sitting in the original fibration $\overline {S^3\setminus {\eta(K)}}={S \times J}/ {h}$
and  $S_1$ is a fiber surface of the fibration.

\begin{lem} \label{lem:converse} Let $M':=\overline {S^3\setminus {\eta(K')}}=S'\times J/ h'$ be an oriented fibered knot complement.
Let $(\Sigma,\   g')$ denote the $HN$-model associated to the fibration with $Q$
the corresponding $HN$-surface of $M'$ sitting in the fibration so that
$S'_1:=S'\times \{{1\over 2}\}$ is disjoint from it.
Let $(\Sigma,\ g'')$ be a second $HN$-model of $M'$ and let
$Q'$
denote the corresponding  $HN$-surface.
Suppose that there  exists an orientation-preserving homeomorphism
$F: M' \longrightarrow M'$ with, $F|\partial M'={\rm id}$, such that
$$ F(Q)=Q' \ {\rm and} \ F(S' \times {x})=S' \times {x}, \  {\rm for\  all}\  x \in J.\eqno (4.6)$$
Then,
there is an orientation-preserving homeomorphism $f: \Sigma \longrightarrow \Sigma$
such that in the mapping class group $\Gamma=\Gamma(\Sigma)$ we have
$$g''=fg' f^{-1}. \eqno(4.7)$$

 \end{lem}

\proof
The existence of the homeomorphism  in (4.6) implies 
that $Q'$ is the $HN$-surface corresponding
to a fibration of $M'$ with fiber 
$S'_1$.  We will now discuss a model of this fibration: If we let $f_1$ denote the restriction of $F$ on the fiber $S'_1$ then the monodromy of our second fibration should be a conjugate of $h'$ by $f_1$ (Proposition \ref{pro:510}). That is the monodromy
of the fibration in which 
$Q'$ is the corresponding $HN$-surface is $h_1:=f_1h'f_1^{-1}$ (where, recall, the equality is understood up to isotopy on the fiber.)
Following the process described in  (4.1)-(4.4)  we can identify $M'$ with
$H\cup_{g_1} (-N) $ where $(\Sigma,\   g_1)$ is the model corresponding to the fibration with monodromy $h_1$.
By Lemma \ref{lem:model510} there is an orientation preserving homeomorphism $f: \Sigma \longrightarrow \Sigma$
such that in the mapping class group $\Gamma=\Gamma(\Sigma)$ we have
$$g_1=fg f^{-1}. $$
Now $M'=H\cup_{g_1} (-N)= H\cup_{g''} (-N)$, with $Q'=\Sigma\cup_{g_1}\Sigma=\Sigma \cup_{g''}\Sigma$ being the $HN$-surface in both splittings. This defines
a homeomorphism  $m:  H\cup_{g_1} (-N) \longrightarrow H\cup_{g''} (-N) $ with $m|\partial M={\rm id}$ and $m|S'_1={\rm id}$  and $m(Q')=Q'$.
Let $m_1$ denote the restriction of $m$ on the  $\Sigma\subset \partial H$ and let $m_2$ denote the restriction of $m$ on $\Sigma=\partial N$. Clearly we have
$g''=m_1^{-1} g_1 m_2$.  
Let $R:=M' \setminus S'_1\cong S'_1\times J$.  
Now $m$ gives rise to a homeomorphism $m: R\longrightarrow R$ such that: (i)  $m(S'_1 \times {x})=S' _1\times {x}$,  for all  $x \in J$; (ii) $m| \partial R={\rm id}$; 
and  (iii) $m(Q')=Q'$. Now  $m$ can be isotopied to  the identity on $R$ by an isotopy that is level preserving (  Lemma 3.5, \cite{kn:Waldhausen}).
Such an isotopy will preserve $Q'$.
It follows that $m_1, m_2$ are isotopic to the identity on $\Sigma$. Since, as discussed earlier,  $g''={m_1}^{{-1}} g_1 m_2$,
$g''=g_1=fg f^{-1}$ up to isotopy in $\Sigma$. \qed
\smallskip

\section{Crossing changes and Dehn twists}
In this section  will prove Theorem \ref{theo:nugatory}. In fact we will work in a more general context as we will consider ``generalized crossing changes".

\subsection{Nugatory crossing changes in fibered knots}
Let $K$ be a knot in ${\bf S}^3$ and let $q\in \Z$.
A generalized crossing of order $q$
on a projection of $K$ is a set $C$
of $|q|$ twist crossings
on two strings that inherit
opposite orientations from any orientation of $K$.
If $K'$ is obtained from $K$ by changing
all the crossings in $C$ simultaneously, we will say that
$K'$ is obtained from $K$ by a generalized crossing change of order $q$ (see Figure 1).\
\begin{figure}[ht]
\centerline{\psfig{figure=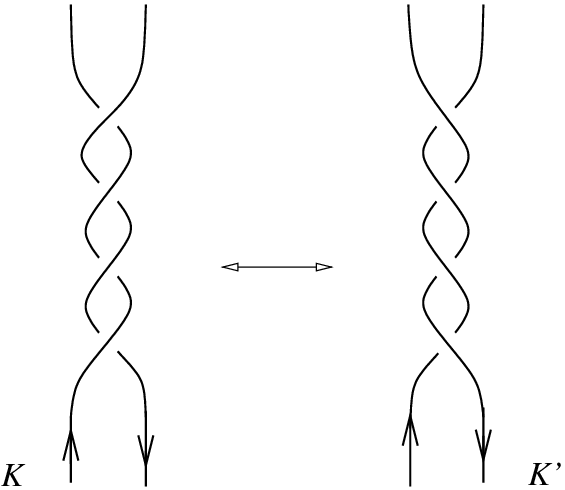,height=1.6in, clip=}}
\caption{The knots $K$ and $K'$ differ by a generalized crossing change
of order $q=-4$.}
\end{figure}
Notice that if $|q|=1$, $K$ and $K'$ differ by an ordinary crossing change
while if $q=0$ we have $K=K'$.
A crossing disc for $K$
is an embedded disc $D\subset S^3$
such
that $K$ intersects ${\rm int}(D)$ twice with
zero algebraic intersection number.
Performing ${\textstyle {{1\over {-q}}}}$-surgery on $L:=\partial D$, for $q\in \Z$,
changes $K$ to another knot $K^{'}\subset S^3$.
Clearly
$K^{'}$ is obtained from $K$
by a generalized crossing change of order $q$. The boundary $L:=\partial D$ is called a
crossing circle supporting the generalized crossing change.

\begin{defi} \label{defi:nugat} A generalized crossing  supported on a crossing circle $L$
of a knot $K$ is called nugatory if and only if
$L:=\partial D$ bounds an embedded disc in the complement of $K$. This disc and $D$ form an
embedded
2-sphere that decomposes $K$ into a connected sum
where some of the summands may be trivial.
\end{defi}

Clearly, changing a nugatory crossing doesn't
change the isotopy class of a knot.
It is an open question whether, in general,
the converse is true (Problem 1.58, \cite{kn:kirby}).
The answer is known to be {\em yes} in the case when $K$ is the unknot
\cite{kn:st} and when $K$ is a
2-bridge knot \cite{kn:to}.
To these we add the 
following theorem.

\begin{theo} \label{theo:fibered'}
Let $K$ be a fibered knot and let $K'$ a
knot obtained from $K$ by a generalized crossing change. If $K'$ is
isotopic to $K$ then a crossing circle $L$
supporting this crossing change bounds an embedded disc in
the complement of $K$.
\end{theo}

\subsection{Preliminaries}
Let $C$ be a generalized crossing of order $q\neq 0$ of a fibered knot $K$. Let
$K'$ denote the knot obtained from
$K$ by changing $C$ and let
$D$ be a
crossing disc for  $C$ with $L:=\partial D$. 

\begin{lem} \label{lem:irreducible} Suppose that  $M_L:=\overline {S^3\setminus {(\eta(K)\cup \eta(L))}}$ is reducible. Then $L$ bounds a disc in the complement of $K$. Thus, in particular,
the crossing change from $K$ to $K'$ is nugatory.
\end{lem}
\proof Let $\Delta$ be an essential 2-sphere in $M_L$; $\eta(K)$ and $\eta(L)$ must lie in different components of $M_L\setminus \Delta$. Isotope $\Delta$ so that its intersection with $D$ is minimal in $M_L$. Then $\Delta \cap D$ is a collection of simple closed curves, each parallel to $\partial D$ on $D$. Thus $K \cup L \subset S^3$ is a split link. Since $L$ is unknotted, it bounds a disc in the complement of $K$.
 \qed
\vskip 0.07in

In the view of Lemma \ref{lem:irreducible}, we may assume that $M_L$ is irreducible.
Since the linking number of $L$ and $K$ is zero, $K$ is homologically trivial in the complement of $L$. It is known that this implies that $K$ bounds a Seifert surface in the complement of $L$. Let
$S$ be a Seifert surface that is of minimum genus among all such Seifert surfaces. Since $S$ is incompressible, after an isotopy we can arrange so that the closed components of  $S\cap D$ are homotopically essential in $D\setminus K$.
But then each such component is parallel to $\partial D$ on $D$  and by further modification we can arrange so that
$S\cap D$ is an arc that is properly embedded on $S$.    
The surface
$S$
gives rise
to Seifert surfaces $S$ and $S'$ of $K$ and $K'$, respectively.

\begin{pro} \label{pro:minimum} Suppose that $K$ is isotopic to $K'$. Then,
$S$ and $S'$ are Seifert surfaces of minimal genus for $K$ and $K'$, respectively.
\end{pro}
\proof We can consider the surface $S$ properly embedded in $M_L$
so that it is disjoint from the component
$\partial \eta(L)$ of $\partial M$. The assumptions
on irreducibility of $M_L$ and 
on the genus of $S$ imply that the foliation machinery of Gabai \cite{kn:ga} applies. In particular,
$S$ is taut in the Thurston norm. The manifolds
$M:= \overline {S^3\setminus {\eta(K)}}$ and $M':=\overline {S^3\setminus {\eta(K')}}$ are obtained by Dehn filling of $M_L$ along $\partial (\eta(L))$.
By Corollary 2.4 of \cite{kn:ga}, $S$ can fail to remain taut in the Thurston norm
(i.e. genus minimizing)
in at most one of $M$ and $M'$. But since $K$ is isotopic to $K'$, $M$ is homeomorphic to $M'$. Thus $S$  
remains taut in both of $M$ and $M'$. This implies that $S$
and $S'$ are Seifert surfaces of minimal genus for $K$ and $K'$, respectively.
\qed

\vskip 0.04in

Next we restrict to fibered knots and recall the assumptions that we have to work with
from the statement of Theorem \ref{theo:fibered'}:
$K$ and $K'$ are fibered knots that are isotopic.
$S$ and $S'$  
are minimum genus Seifert
surfaces,  for $K$ and $K'$, respectively.

%By the discussion above, $M'=\overline {S^3\setminus {\eta(K')}}$ is obtained %from $M=\overline {S^3\setminus {\eta(K)}}$ by performing ${\textstyle {{1\over %{-q}}}}$-surgery on $L:=\partial D$.

\smallskip

\subsection{An $HN$-model for $M'$ from Dehn surgery.} 
With the notation of Section 4, 
there is a fibration
$M:=\overline {S^3\setminus {\eta(K)}}= {{S\times J}/ {h}}$ 
with monodromy 
$h: S \longrightarrow S$.  With $N:=S\times [0, 1]$
and $\Sigma:=\partial N$ we have an
$HN$-model $(\Sigma, g)$ corresponding to the fibration of $M$.
We can think of the Heegaard splitting of $S^3$ corresponding to the fibration $M={{S\times J}/ {h}}$ 
as the quotient

$$N\cup_g (-N):= N{ \sqcup} (-N) / \{ y\sim ig(y) |  y\in \Sigma \}. \eqno(5.1)$$
We will further assume that the crossing circle $L$ is embedded
on $\Sigma$ so that $D$ is a meridian disc of $N$. We will furthermore assume that
the embedding of $L$ on $\Sigma$ is chosen so that, up to isotopy
in $M$, the geometric intersection $|K\cap D|$ is minimal. Note 
that since we assumed that
$M_L:=\overline {S^3\setminus {(\eta(K)\cup \eta(L))}}$ is irreducible, this minimum intersection must be non-zero.
Let $\tau: N \longrightarrow N$ denote the right-handed Dehn twist of $N$ along the meridional disc $D$
and let $T_L:=\tau| \Sigma$, where $L=\partial D$. We have
$\tau^{-q}(S)=S'$ and $\tau^{-q}(K)=K'$.
Recall that 
$M:={{S\times J}/ {h}}$ and that
 $M':=\overline {S^3\setminus {\eta(K')}}$ is obtained from $M_L$ by
Dehn filling along $\partial{\eta(L)}$ with slope 
${\textstyle {{ {{1\over {-q}}}}}}$.  Next we use that information to construct
an $HN$-model for $M'$. The proof of Lemma \ref{lem:model}
follows a known process of passing between gluing maps of Heegaard splittings and Dehn surgeries of 3-manifolds
(compare,  pp. 86-87 of \cite{kn:amc}).

\begin{lem} \label{lem:model} $(\Sigma, \ gT^{-q}_L)$ is an  $HN$-model for $M'=\overline {S^3\setminus {\eta(K')}}$.
\end{lem}

\proof By assumption 
$(\Sigma, g)$ is an $HN$-model corresponding to the fibration $M= {{S\times J}/ {h}}$.
Let $A$ denote an annulus on $\Sigma$ that
supports $T_L$ and let $B:=g(A)$.
We will think of this $HN$-splitting of $M$ as the quotient

$$H\cup_g (-N):= H{ \sqcup} (-N) / \{ y\sim ig(y) |  y\in \Sigma \}, \eqno(5.2)$$
where $H\subset N$.
We consider the complement $M_L:=\overline{S^3\setminus (\eta(K)\cup \eta(L))}$ as the pre-quotient space
$$H \cup_{g^{1}} (-N) \ \ \ {\rm where} \ \ \ g^1:= g| (\Sigma \setminus A):
\Sigma \setminus A \longrightarrow \Sigma \setminus B. \eqno(5.3) $$
Thus  we can think of the torus ${\mathcal T}:=A\cup B$ as the boundary
torus of a tubular neighborhood of $L$. Let $\alpha$ be an arc that is properly embedded and essential on $A$ such that it intersects $L$ exactly once and let $\beta:=g(\alpha)$.
Now $\mu:=\alpha \cup \beta$ is the meridian of ${\mathcal T}$ and $\lambda:=L$ is the longitude which we will
orient so that their algebraic intersection number on ${\mathcal T}$, denoted by $<\lambda, \mu>$, is one.
Since $T_L$ is supported in $A$ it can be considered as a Dehn twist on ${\mathcal T}$.
We have 
$$T_L^{q}(\mu)=\mu-q\lambda=T_L^{q}(\alpha)\cup \beta.$$
(Recall that, in general, if $a,b$ are simple closed curves on ${\mathcal T}$, we have
$T_a(b)=b+ <a,b> a$. Since $<\lambda, \mu>=1$, we have $T_L^{-1}(\mu)=\mu+\lambda$, which explains the change of sign
between the power $T_L^q$  and the coefficient of $\lambda$ in $T_L^q(\mu)$ in the equations above.)

Now if we  have
$$\alpha'\cup \beta= \mu-q\lambda.$$
Let $M_L(q)$ denote the 3-manifold obtained from $M_L$
by ${1}\over {-q}$ \ Dehn filling on $\mathcal T$.
From the discussion above, in order to obtain $M_L(q)$ one needs to attach a
solid torus to $\mathcal T$ in such a way so that the meridian is attached along the
curve  $\mu$. It follows that
$H\cup_ {g T_L^{-q}} (-N)$ is an $HN$-splitting for $M_L(q)$.
But since by assumption we have
$M_L(q)=\overline {S^3\setminus {\eta(K')}}=M'$, it follows that $(\Sigma, \ gT^{-q}_L)$ is an $HN$-model for $M'$. \qed
\vskip 0.04in

\subsection{ Understanding the $HN$-model $(\Sigma, \ g T_L^{-q})$} In the view of
the conventions adapted earlier, $N$ is thought as a product $S\times I$
and
$K$ is embedded on $\Sigma:=\partial N$. The Dehn twist 
$\tau^{-q}: N \longrightarrow N$
changes $K$ to $K'$ and the product structure of $N$ to $S'\times I$.
By our assumptions, each of $K$, $K'$ split $\Sigma$ into two bounded surfaces that
are incompressible in $N$.
Let $A$ be an annulus on $\Sigma$ supporting the restriction
$T_L:= \tau | \Sigma$ so that the core of $A$ is $L$ and the intersection $A\cap K$ consists of two  properly embedded, disjoint arcs, say $\alpha_{1}, \alpha_{2}$,
each of which intersects $L$ exactly once.
We set, $B:=g(A)$,
$\gamma:=g(K)$, $\gamma':= g(K')$ and $z:=g(L)$. By construction, we have $g|K={\rm id}$. Thus, $g^{-1}(K)=K$, $B\cap \gamma =\alpha_1 \cup \alpha_2$.
We have $$\gamma':=g(K')=g(T^{-q}_L(K))= g(T^{-q}_L ( g^{-1}(K)))
=gT^{-q}_L g^{-1}(K)=T^{-q}_{g(L)}(K),$$
where the last equation follows from the fact that in the mapping class group we have  
$gT_L g^{-1}=T_{g(L)}$.
Thus $\gamma'$ is the result of $\gamma:=g(K)=K$ under a non-trivial power of a Dehn twist along $z:=g(L)$ supported on $B$.  
We will think of the $HN$-splitting of
$M'=\overline {S^3\setminus {\eta(K')}}$ corresponding to the model $(\Sigma, \ gT_L^{-q})$ as the quotient
$$M'= H{ \cup} (-N) / \{ y\sim igT^{-q}_L(y)|   y\in \Sigma \}, \eqno(5.4)$$
\noindent and we will identify the corresponding Heegaard
splitting of $S^3$ with the quotient
$$N{ \sqcup} (-N) / \{ y\sim igT^{-q}_L(y) |   y\in \Sigma \}. \eqno(5.5)$$
\vskip 0.04in

%Next we observe that the irreducibility assumption implies that the Dehn twist 
%$T_{g(L)}$ does not extend over the handlebody $N$.

%\begin{lem}\label{noextention} Suppose that there exists a homeomorphism $\Phi: N \longrightarrow N$
%such that $\Phi | \Sigma=T_{g(L)}$. Then, $g(L)$ bounds an embedded disc in $N$ and $M_L$
%is reducible.
%\end{lem}
%\proof
% Theorem
%\ref{theo:mac} implies that $g(L)$ must bound a disc, say $D'$, in $N$.
%Recall that in the Heegaard splitting of (5.1)
 %$g(L)$ and $L$ are identified and that $L$ bounds a disc $D$ in $N$.
%Now the 2-sphere
 %$D\cup D'$ intersects the Heegaard surface of (5.1) at $L$ and the knot $K$ at exactly 2-points. This sphere realizes $K$ as a (possibly trivial)
%connect sum and it contains $L$. It follows that $L$ bounds an embedded disc in the complement of
%$K$ and $M_L$ reducible. \qed

\begin{lem} \label{lem:greducible} Push $g(L)$ slightly in the interior of $N$
an let $\eta(g(L))$ denote a tubular neighborhood of it in there.
If $\overline{N\setminus \eta(g(L))}$ is reducible,
then $M_L:=\overline{S^3\setminus{\eta(K)\cup \eta(L)}}$ is reducible.

\end{lem}

\proof Since $\overline{N\setminus{ \eta(g(L))}}$
 is reducible, $g(L)$ must lie in a 3-ball in $N$. It follows that $K \cup L \subset S^3$ is a split link, thus $M_L$ is reducible.\qed
\smallskip

\vskip 0.06in

In the view of Lemma \ref{lem:greducible} and our earlier assumption that $M_L$ is
irreducible we may assume that $\overline{N\setminus{ \eta(g(L))}}$ is irreducible.

For $i=0,1$, let $S_i:=S\times \{i\}$.
The boundary $\partial N$  is the union $S_0\cup E \cup S_1$,
where $E=\partial S \times (0, \ 1)$. Let $\Sigma_0, \Sigma_1$
denote the image of $S_0, S_1$, respectively, 
under the Dehn twist
$T^{-q}_{g(L)}$. Then, for $i=0,1$, $\partial \Sigma_i= \gamma'\times \{i\}$.
\smallskip

\begin{lem} \label{caseA} The surfaces $\Sigma_0, \Sigma_1$ are incompressible in $N$.
\end{lem}
\proof Suppose, on the contrary that one of  $\Sigma_0, \Sigma_1$, say
 $\Sigma_0$
compresses in $N$.  Consider $N$ as a product $S\times I$ with $g(L)$ a knot in $N$.
By assumption $\Sigma_0$ compresses in $N$. Performing the Dehn twist
$T^{-q}_{g(L)}$ is equivalent to doing surgery along $g(L)$. Since $q\neq 0$, this surgery is non-trivial (Proposition  \ref{pro:nevertrivial}). Now $\Sigma_0$ is the result of $S_0$
under this surgery. Thus there is a non-trivial surgery in $S\times I$ such that the surface $S_0$ compresses in the manifold obtained after surgery.
By Theorem 1 of \cite{kn:st1} there is a simple closed homotopically essential
curve $L' \subset \Sigma_0$ such that $g(L)$ and $L'$ cobound an embedded annulus in  $N=S\times I$.
Furthermore, this annulus determines the slope of the surgery. This implies that $g(L)$ bounds a disc in $N$. But then,
any Dehn twist on $\partial N$ along $g(L)$ extends to a Dehn twist on $N$; a homeomorphism of $N$.
Since $S_0, S_1$ are incompressible, their images under any homeomorphism of $N$ are also incompressible in $N$. This contradicts the assumption that $\Sigma_0$ compresses. \qed
\smallskip

\begin{lem} \label{lem:structure} With the notation and the setting as above, there exists a fibration of $M'$, with fiber $S'$
and corresponding $HN$-model $(\Sigma, g_1)$, and an orientation preserving 
homeomorphism $f: \Sigma \longrightarrow \Sigma$
such that in $\Gamma(\Sigma)$ we have
$$g'':=gT^{-q}_L=fg_1 f^{-1}. \eqno(5.6)$$
\end{lem}
\proof 
We recall that the Heegaard splitting in (5.5) is the result of the splitting of (5.1)
after the Dehn twist $\tau^{-q}$ on $N$. This twist changes
the product structure of $N$ from $S\times I$ to $S' \times I$.  For $i=0,1$, let $S_i:=S\times \{i\}$.
The boundary $\partial N$  is the union $S_0\cup E \cup S_1$,
where $E=\partial S \times (0, \ 1)$.  We have
$$g''(S_i)=g(T^{-q}_L(S_i))= g(T^{-q}_L ( g^{-1}(S_i)))
=gT^{-q}_L g^{-1}(S_i)=T^{-q}_{g(L)}(S_i)=\Sigma_i.$$
By Lemma   \ref{caseA}, $\Sigma_i$ is incompressible in $N$.
Now we pass to the corresponding $HN$-splitting of
(5.4) and we 
use $Q'$ to denote the corresponding $HN$-surface.
Since the $HN$-surface
of $H\cup_g (-N)$ is disjoint from a level surface  of the fibration
$S\times J/h$,
$Q'$ is disjoint from a neighborhood of a copy $S' \subset {\rm int}(H)$.
By Proposition \ref{pro:510}, $M'$ is fibered with fiber $S'$.
Let  $(\Sigma, g_1)$ denote the $HN$-model of this fibration and let $Q$ denote the corresponding 
$HN$-surface.  On one hand  $M'$ cut along $S'$ is a product $S'\times J$.
On the other hand  $M'=\overline{S^3\setminus \eta(K)}$ is homeomorphic to
$$H{ \sqcup} (-N) / \{ y\sim ig''(y) |   y\in \Sigma \}, \eqno(5.7)$$ 
For  $i=0,1$, the surface $S_i\cup_{g''}\Sigma_i \subset Q'$
gives a properly embedded incompressible surfaces in $M'$. These two surfaces  can be isotopied in $M'$, relatively $\partial M'$,  so that  each becomes parallel to
the fiber $S'$ (Proposition 3.1 of \cite{kn:Waldhausen}). In fact, the isotopy brings each of the surfaces onto a level surface of the fibration
(Proposition  \ref{pro:510}).
This implies that there is an orientation preserving homeomorphism $F: M' \longrightarrow M'$ with, $F|\partial M'={\rm id}$, such that
$ F(Q)=Q' \ {\rm and} \ F(S' \times {x})=S' \times {x}, \  {\rm for\  all}\  x \in J$.
Now applying Lemma \ref{lem:converse}
to the models $(\Sigma, g_1)$ and $(\Sigma, g'')$ we get the desired conclusion.
\qed

%Equivalently
%we have a fibration $M'=S'\times J/h_1$, with monodromy $h_1: S' \longrightarrow S'$. We will denote
%by $Q$ the $HN$-surface corresponding to this fibration. The monodromy gives rise to an $HN$-model  $(\Sigma, g_1)$. 
%We will apply Lemma \ref{lem:converse} to 
%$(\Sigma, g_1)$, $(\Sigma, \ g'')$: By Proposition \ref{pro:510}, $S'_1$ can be taken to be a level surface
%of the  fibration $M'=S'\times J/h_1$. Since $S'_1$ was chosen so it is disjoint from $Q,Q'$,
%condition (1) of . 
%Let $R'$ denote $M'$ cut along $S'_1$; we have $R'=S'_1\times J$.

%There is a
%homeomorphism of $M'$ that takes $\Sigma'_1$ to the fiber $S'$; in particular $M'$ cut along 
%$\Sigma'_1$  is a product $\Sigma'_1\times J$ . We can now find a 
%a level preserving homeomorphism $F: M' \longrightarrow M'$ such that 
%and $F$ maps $Q$ to $Q'$ and so that
%$F(S'_1)=S'_1$. Thus Lemma \ref{lem:converse}
%applies to give the desired conclusion.  \qed

\smallskip

\subsection{Proof of Theorem \ref{theo:fibered'}} Let $K,K'$ be fibered isotopic knots, such that $K'$ is obtained from
$K$ by a generalized crossing change, of order $q\neq 0$, supported on a crossing circle $L$.
Let $D$ be a crossing disc with $L:=\partial D$. 
We will consider the Heegaard splittings of (5.1) and (5.5)
so that the crossing circle $L$ is embedded
on $\Sigma$ and $D$ is a meridian disc of $N$. Recall that the crossing change from $K$ to $K'$ is now achieved
by the Dehn twist $\tau^{-q}$ of $N$ along $D$.

We will assume that
$L$ is homotopically essential on $\Sigma$ since otherwise the crossing change from $K$ to $K'$ is obviously nugatory. 

If 
$M_L:=\overline{S^3\setminus{\eta(K)\cup \eta(L)}}$ is reducible, then we are done by Lemma \ref{lem:irreducible}. We will assume that
$M_L:=\overline{S^3\setminus{\eta(K)\cup \eta(L)}}$  is irreducible. Then, by
Lemma \ref{lem:greducible}  $\overline{N\setminus{ \eta(g(L))}}$ is irreducible.
By Proposition \ref{pro:minimum}, $S$ and $S'$ are of minimum genus for $K$ and $K'$, respectively.
By Lemma \ref{lem:structure}, there is an $HN$-model $(\Sigma, g_1)$ that corresponds to a fibration
$M'=S'\times J/h_1$ and $f: \Sigma \longrightarrow \Sigma$ so that
(5.6) is satisfied. Equivalently, we have
$f^{-1}gT^{-q}_Lf=g_1$.
Since $K$ and $K'$ are isotopic knots there is an orientation preserving homeomorphism,
say $\Phi$, of $S^3$
that brings $K$ to $K'$. 
Now we have two equivalent fibered knot complements;  $M'=S'\times J/h_1$ and $M=S\times J/h$.
Via Lemma \ref{lem:model510}, $\Phi$ gives rise to
a homeomorphism $\phi: \Sigma \longrightarrow \Sigma$ such that

$$gT^{-q}_L=\phi g \phi^{-1} \ \ {\rm or} \ \
T^{-q}_L=g^{-1}\phi g \phi^{-1}. \eqno(5.8)$$
Now  (5.8) realizes
$T_L^{-q}$ as a  commutator of length one in $\Gamma$.
By Corollary \ref{corol:one}, $L$ must be homotopically trivial on $\Sigma$ which contradicts the assumption that
$M_L$ is irreducible.
\qed

%\smallskip

Since Kotschick's result is not true in the case of twists with mixed signs, the argument above breaks down in an attempt to generalize the statement of Theorem \ref{theo:fibered'}.
to multiple crossing changes. But as the following example shows the result is, in fact, not true!

\begin{exam} \label{exam:figure8}{\rm Let $K$ denote the
figure eight knot as boundary of  a genus one Seifert surface $S$ obtained by Hopf plumbing two
once twisted bands $B_L$ and $B_R$. Consider $D_1$, $D_2$ crossing discs of $K$ such that
$D_1\cap B_L$ (resp. $D_2\cap B_R$) is an essential arc cutting $B_L$ (resp. $B_R$)
into a square. One can perform opposite sign twists of order four along $D_1$, $D_2$ to
transform $S$ to $S'$ where in $S'$ the Hopf band $B_L$ becomes the Hopf band $B_R$
and vice versa. The knot $K':=\partial S'$ is isotopic to $K$. Moreover,  $S$ and $S'$
are clearly minimum genus
Seifert surfaces for $K$ and $K'$, respectively. However,
neither of $L_1:=\partial D_1$ or $L_2:=\partial D_2$
bounds disc in the complement of $K$.}
\end{exam}

\smallskip

\section{Adjacency to fibered knots}

We begin by recalling from \cite{kn:kl} the following definition.

\begin{defi} \label{defi:adj} Let $K$, $K'$ be knots.
We
will say that $K$ is $n$-adjacent to $K'$, for some
$n\in {\bf N}$,
if $K$ admits a projection containing $n$ generalized crossings
such that changing any $0<m\leq n$ of them yields a
projection of $K'$.
\end{defi}

In \cite{kn:kl} we showed the following: Given knots $K$ and $K'$ there exists a constant $c=c(K,K')$ such that if $K$ is $n$-adjacent to $K'$ then either $n\leq c$ or $K$ is isotopic to $K'$.
Here, using Theorem \ref{theo:fibered'},
we will show that if $K'$ is assumed to be fibered, then we can have a much stronger result.

\begin{theo}\label{theo:fibered} Suppose that $K'$
is a fibered knot and that $K$ is a knot such that $K$ is $n$-adjacent to $K'$, for some $n>1$.
Then, either $K$ is isotopic to $K'$ or we have  $g(K)>g(K')$.
\end{theo}

\begin{remark} \label{remark:forsmall} {\rm It is not hard to see that if $K$ is $n$-adjacent to $K'$,
for some $n>1$, 
then $K$ is $m$-adjacent to $K'$, for all $0< m\leq n$.}
\end{remark}

Suppose that $K$ is $n$-adjacent to $K'$ and let $L$ be a
collection of $n$ crossing circles supporting the
set of generalized crossings that exhibit $K$ as $n$-adjacent to $K'$.
Since the linking number
of
$K$ and every component of $L$ is zero,
$K$ bounds a Seifert surface $S$ in the complement
of $L$. Define
$$g_n^L(K):={\rm min}\,\{\, {\rm genus (S)}\,|\,S\, {\rm a\ Seifert \ surface\  of \,}K\, {\rm as\ above}\,\}.$$
We recall the following.

\begin{theo} \label{theo:genus} {\rm [Theorem 3.1, \cite{kn:kl}]}
We have
$$g_n^L(K)={\rm max}\,\{\, g(K),g(K')\,\}$$
where $g(K)$ and $g(K')$ denote the genera of $K$ and $K'$, respectively.
\end{theo}

\proof  {\rm [Proof of Theorem \ref{theo:fibered}]}
Let $K'$ be a fibered knot. In the view of Remark \ref{remark:forsmall},
it is enough to prove that
if $K$ is a knot that is 2-adjacent to $K'$ then
either $K$ is isotopic to $K'$ or we have $g(K)>g(K')$.
To that end, suppose that $K$ is exhibited as 2-adjacent
to $K'$ by a  two component crossing link $L:=L_1\cup L_2$.
Let $D_1,D_2$ be crossing discs for $L_1,L_2$, respectively.
Suppose, moreover, that $g(K)\leq g(K')$; otherwise
there is nothing to prove.
Let $S$ be a Seifert surface for $K$ that is of minimal
genus among all surfaces bounded by $K$ in
the complement of $L$. As explained earlier in the paper,
we can isotope $S$ so that, for $i=1, 2$,
$S\cap {\rm int}(D_i)$
is an arc, say $\alpha_i$ that is properly embedded in $S$.
For $i=1,2$, let $K_i$ (resp. $S_i$) denote the knot (resp. the Seifert surface)
obtained from $K$ (resp. $S$) by changing $C_i$. Also let $K_{3}$
denote the knot obtained by changing $C_1$ and $C_2$ simultaneously and let $S_{3}$ denote the corresponding surface. By assumption, for $i=1,2,3$,
$K_i$ is isotopic to $K'$ and $S_i$
is a Seifert surface for $K_i$. Since
$g(K)\leq g(K')$,
Theorem
\ref{theo:genus} implies that $S_i$
is a minimum genus surface for $K_i$.
Observe that $K_3$ is obtained from $K_1$ by changing $C_2$ and that
they are fibered isotopic knots.
Furthermore, $S_3$ is obtained from $S_1$
by twisting along the arc $\alpha_2 \subset
S$. By Theorem \ref{theo:fibered'}, $L_2$ bounds an embedded disc
$\Delta_2$ in the complement of $K_1$. Since $S_3$ is incompressible, after an isotopy,
 we can assume that
$\Delta \cap S_3=\emptyset$. Now let us consider the 2-sphere
$$S^2:= \Delta \cup D_2.$$ By assumption $S^2 \cap S_3$
consists of the arc $\alpha_2 \subset S_3$. Since $\alpha_1$ and $\alpha_2$
are disjoint, the arc $\alpha_1$ is disjoint from $S^2$.
But since $K$ is obtained from $K_1$ by twisting along $\alpha_1$,
the circle $L_2$ still bounds an embedded disc in the complement of $K$.
Hence, $K$ is  isotopic to $K'$. \qed

\begin{remark}\label{remark:noteq} {\rm
The trefoil knot is 2-adjacent to the unknot. Since the trefoil is a fibered knot Theorem \ref{theo:fibered} implies that
the unknot is not 2-adjacent to the trefoil.
Thus $n$-adjacency is not an equivalence relation on the set of knots.}
\end{remark}

\end{document}